\documentclass[a4paper, 11pt]{amsart}
\usepackage[T1]{fontenc}
\usepackage[latin1]{inputenc}
\usepackage[english]{babel}
\usepackage{amsmath}
\usepackage{amssymb}
\usepackage{amsfonts}

\begin{document}
\newcommand{\normi}[1]{\|#1\|}
\newcommand{\itse}[1]{\left|\,#1\right|}
\newcommand{\its}[1]{\bigl|\,#1\bigr|}
\newcommand{\rn}{\mathbb{R}^n}
\newcommand{\na}{\mathbb{N}}
\newcommand{\re}{\mathbb{R}}
\newcommand{\R}{\mathcal{R}}
\newcommand{\Li}{\mathcal{L}}
\newcommand{\Z}{\mathbb{Z}}
\newcommand{\M}{\mathcal{M}}
\newcommand{\vM}{\overset{\rightarrow}{M}}
\newcommand{\ve}[1]{\overset{\rightarrow}{#1}}
\def\Xint#1{\mathchoice
 {\XXint\displaystyle\textstyle{#1}}%
 {\XXint\textstyle\scriptstyle{#1}}%
 {\XXint\scriptstyle\scriptscriptstyle{#1}}%
 {\XXint\scriptscriptstyle\scriptscriptstyle{#1}}%
 \!\int}
 \def\XXint#1#2#3{{\setbox0=\hbox{$#1{#2#3}{\int}$}
 \vcenter{\hbox{$#2#3$}}\kern-.5\wd0}}
 \def\ddashint{\Xint=}
 \def\dashint{\Xint-}
\newcommand{\m}{m}
\newcommand{\subsub}{\subset\subset}
\newcommand{\av}[3]{\underset{B(#2,#3)}{\dashint}#1(y)\,dy\,}
\newcommand{\avr}[1]{\underset{#1}{\dashint}}
\newtheorem{theorem}{Theorem}[section]
\newtheorem{lemma}[theorem]{Lemma}
\newtheorem{lause}[theorem]{Theorem}
\newtheorem{proposition}[theorem]{Proposition}
\newtheorem{definition}[theorem]{Definition}
\newtheorem{corollary}[theorem]{Corollary}
\newtheorem{question}[theorem]{Question}
\newtheorem{remark}[theorem]{Remark}
\title[Differentiability of directionally differentiable functions]
{On the differentiability of directionally differentiable functions and applications}
\author{Hannes Luiro}
\address{Department of Mathematics and Statistics\\
University of Jyvï¿½skylï¿½\\
P.O.Box 35 (MaD)\\
40014 University of Jyvï¿½skylï¿½, Finland}
\email{haluiro@maths.jyu.fi}
\maketitle
\begin{abstract}
In the first part of this paper we establish, in terms of so called $k$-tangential sets, a kind of optimal estimate for the size and structure of the set of non-differentiability of Lipshitz functions with one-sided directional derivatives. These results can be applied to many important special functions in analysis, like distance functions or different maximal functions.

In the second part, having the results from the first part in our use, we focus more carefully on the differentiability properties of the classical Hardy-Littlewood maximal function. For example, we will show that if $f$ is continuous and differentiable outside a $\sigma$-tangential set, then the same holds to the maximal function $Mf$ as well (if $Mf\not\equiv\infty$). As an another example, our results also imply that if $f$ is differentiable almost everywhere (and $Mf\not\equiv\infty$), then $Mf$ is differentiable a.e. 
\end{abstract}

\section{Introduction}
The original motivation of this research is raised by the study of the differentiability properties of the classical Hardy-Littlewood maximal function
\begin{equation}\label{kaava:1}
\M f(x)=\sup_{r>0}\avr{B(x,r)}|f(y)|\,dy\,=\sup_{r>0}
\frac{1}{|B(x,r)|}\int_{B(x,r)}|f(y)|\,dy\,.
\end{equation}
The first step on this field of study was taken by J. Kinnunen who observed that $M$ is a bounded operator in the Sobolev-spaces $W^{1,p}(\rn)$ when $1<p\leq\infty$ \cite{Ki}. Some of the extensions and related results can be found e.g. from \cite{AP},\cite{HM},\cite{HO},\cite{KL}, \cite{Ko}, \cite{KS}, \cite{Lu} and \cite{Lu2}. 

Our initial goal, in the case of maximal functions, was to find optimal estimates for the size of the non-differentiability points of $Mf$ in the case where $f$ itself is differentiable or a.e. differentiable. It seems that in the context of maximal operators these questions have not been studied before. Maybe the most relevant result in this direction is from P. Haj{\l}asz and J. Mal\'{y} who showed that $Mf$ is \textit{approximately} differentiable if $f$ is approximately differentiable (\cite[Thm. 1]{HM}).

The investigation of the above question quickly led us to consider much more general problem, the differentiability of functions $f:\rn\to\re$ such that the 'one-sided' directional limit
\begin{equation}
\lim_{\lambda\to 0^+}\frac{f(x+\lambda\theta)-f(x)}{\lambda}=:D_{\theta}f(x)
\end{equation}
exists for every $x\in\rn$ and unit vector $\theta\in\rn$. We will call these functions as \textit{directionally differentiable} functions. Many important special functions in analysis, while not being everywhere differentiable, like convex functions, distance functions, pointwise maximum functions, different maximal functions etc., are still directionally differentiable. 
As we will see, this applies also to the Hardy-Littlewood maximal function $Mf$, if $f$ itself is differentiable. 

The interplay between the directional differentiability and the real differentiability, also in more general context of Banach spaces, have been most intensively studied in optimization theory but also in geometric measure theory and PDE-theory (see e.g. \cite{ACP}, \cite{C}, \cite{BC} and references therein). Naturally, the main focus has been on the evaluation of the size and structure of the set of non-differentiability.

It is not difficult to see (see e.g. \cite{ACP}) that the set of non-differentiability of a directionally differentiable Lipschitz-function $f$ is $\sigma$-porous, thus it can be included in a countable union of sets $E_i$ such that for any $x\in E_i$ there exists $0<\delta_i(x)<\frac{1}{2}$ so that for any $0<r<1$ there exists a ball $B(y,\delta r)\subset B(x,r)\setminus E_i$. More careful analysis shows that actually $\delta(x)$ can be chosen independently on $x$ or $i$ and also arbitrarily close to $\frac{1}{2}$. This argument also implies the obviously best Hausdorff dimension estimate, $n-1$, for the set of non-differentiability. 

However, from topological point of view, the above porosity results can be evidently improved, but it seems that more precise estimates have not been published before, except in the most simple (and important) special case of convex functions (see \cite{Z} and \cite{P}). It seems that the main focus in the previous researches on this area (like in \cite{BC}) has not been on proving optimal results in the case of real-valued functions on euclidean spaces but on more general context of Banach spaces.

In our main Theorem \ref{OK} we prove, in terms of $k$-tangential sets, a kind of optimal estimate for the size and structure of non-differentiability points of directionally differentiable Lipschitz functions. Before stating the results, let us introduce some notation. For a linear subspace $V\subset\rn$, its orthogonal complement is denoted by $V^{\perp}$. We say that $E\subset\rn$ is $k$-\textit{tangential}, $0\leq k\leq n-1$, if for every $x\in E$ there exists a $k$-dimensional linear subspace $V_x$ so that if $x+h_k\in E$ and $h_k\to 0$ as $k\to\infty$, then 
\begin{equation*}
\frac{|h_k^{V_x^{\perp}}|}{|h_k^{V_x}|}\to 0\,\text{ as }k\to\infty\,,
\end{equation*}
where $h_k=h_k^{V_x}+h_k^{V_x^{\perp}}$ such that $h_k^{V_x}\in V_x$ and $h_k^{V_x^{\perp}}\in V_x^{\perp}\,$.
We call $n-1$-tangential sets simply as tangential sets. If $E$ can be covered by a countable union of $k$-tangential sets, we say that $G$ is $\sigma$-$k$-\textit{tangential}.  

It is clear that $k$-tangential sets does not need to have finite Hausdorff $k$-measure. However, it is well known that $k$-tangentiality implies the $k$-rectifiability (see e.g. \cite[Lemma 15.13]{M}), which again implies that every $\sigma$-$k$-tangential set can be covered by a countable union of $k$-tangential sets with finite Hausdorff $k$-measure.

In our main Theorem \ref{OK} we will show that outside a $\sigma$-$k$-tangential set, there exists a $k+1$-dimensional linear subspace $V_x$ so that $f$ is differentiable at $x$ in respect to $V_x$ and all 'halfspaces' $H$ of the form \[H=\{v+\lambda b:v\in V_x\,,\,\lambda \geq 0\}\,\,,\, b\in\rn\,.\] The differentiability in respect to certain linear subspace (or halfspace) $V$ here simply means that the restriction of $f$ to $x+V$ is differentiable at $x$. The proof and exact formulation of this result are given in Section 2.  

The following theorem is a direct consequence of Theorem \ref{OK}:
\begin{theorem}\label{ekapeli}
Let $f:\rn\to\re$ be directionally differentiable and Lipschitz. Then $f$ is differentiable up to a countable union of tangential sets.
\end{theorem}

In the end of the Section \ref{general} some basic examples of applications of Theorem \ref{OK} are given.
In Section \ref{last} we go back into our original question concerning the Hardy-Littlewood maximal functions. We prove the following theorem: 
\begin{theorem}\label{differentioituva}
If $f$ is continuous and differentiable outside a $\sigma$-tangential set and $Mf(x)\not\equiv\infty$, then $Mf$ is continuous and differentiable up to a $\sigma$-tangential set.
\end{theorem}
A big part of the proof of this theorem follows from Theorem \ref{OK}, but some additional results are still needed, essentially because we do not assume that $f$ is $C^1$-function nor anything about the behaviour of $f$ at infinity. Some of those auxiliary lemmas have also their own interest. For example, as a corollary of our results we get that maximal operator preserves the almost everywhere differentiability, providing a natural counterpart for \cite[Thm. 1]{HM}.\\

\textit{Acknowledgements.}
The author would like to thank Lizaveta Ihnatsyeva, Juha Kinnunen and Antti Käenmäki for useful conversations.

\section{Directionally differentiable Lipschitz-functions}\label{general} 
\subsection{Premilinaries}
Let us begin with listing some notation. The unit sphere $\{x\in\rn\,:\,|x|=1\}$ is denoted by $S^{n-1}\,$ and unit ball by $B_n$.
For $A\subset\rn$ and $x\in\rn$, define $|x-A|=\inf\{|x-y|:y\in A\}$.
If $\delta>0$, then the $\delta$-extension of $A\subset\rn$ is defined by
\begin{equation*}
A_{(\delta)}=\{x\in\rn:\,|x-A|\leq\delta\,\}.
\end{equation*}
The \textit{Hausdorff distance} $\mathcal{H}$ for $A\subset\rn$ and $B\subset\rn$ is defined by
\begin{equation*}
\mathcal{H}(A,B)=\inf\{\delta>0\,:\,A\subset B_{(\delta)}\text{ and }B\subset A_{(\delta)}\,\}.
\end{equation*}

Moreover, let $A,B\subset\rn$ and define  
\begin{align*}
&A+B=\{a+b\,:\,a\in A, b\in B\,\}\,,\\
&\langle A\rangle=\{\,\sum_{i=1}^k \lambda_ka_k\,:\,k\in\na, a_k\in A\,,\lambda_k\in\re\,\}\,,\\
&\langle A\rangle^+=\{\,\sum_{i=1}^k \lambda_ka_k\,:\,k\in\na, a_k\in A\,,\lambda_k\geq 0\}.
\end{align*}
In the proof of Theorem \ref{OK} we have to do some elementary linear algebra on subspaces of the form $W=\langle A\rangle^+$, called as \textit{semi-linear} subspaces in future. The set of these spaces is denoted by $\mathcal{W}_n$. The subset of $\mathcal{W}_n$ of all linear subspaces is denoted by $\mathcal{V}_n$.

\textit{Remark.} If considered just for the purposes of the results below, the above definitions (as well as some auxiliary lemmas below) are somewhat too general kind of. Indeed, in the proof of the main theorem the most complicated semi-linear subspaces are of the form $H-\langle a\rangle^+$, where $a\in\rn$ and $H$ is a \textit{half-space} of the form $H=V+\langle b\rangle^+$ such that $V\in \mathcal{V}_n$ and $b\in\rn\,$. However, that kind of general approach was chosen, partly because that choice do not remarkably complicate the presentation, partly for the purposes of possible further studies e.g. in finite dimensional Banach spaces.\\ 

We endow $\mathcal{W}_n$ and $\mathcal{V}_n$ with metric $\mathcal{H}_c$ which is just a restriction of the Hausdorff metric to the unit ball, thus  
\begin{equation}
\mathcal{H}_c(A,E)=\mathcal{H}(A\cap B_n,E\cap B_n)\,.
\end{equation}
Observe that $\mathcal{H}_c(A,E)=0$ if and only if the closures $\overline{A\cap B_n}$ and $\overline{E\cap B_n}$ coincide and, moreover, $\mathcal{H}_c(A,B)=\mathcal{H}_c(A',B')$ if $\overline{A}=\overline{A'}$ and $\overline{B}=\overline{B'}\,$. For this reason, we will assume in future, without changing any notation that two elements of $\mathcal{W}_n$ coincide if their closures coincide, thus $\mathcal{W}_n$ is actually the set of equivalence classes determined by that equivalence relation.

Now it is easy to check that $(\mathcal{W}_n,\mathcal{H}_c)$ and $(\mathcal{V}_n,\mathcal{H}_c)$ are compact\footnote{It is well known that for any compact set $K\subset\rn$, the set of all subsets of $K$ is compact, if endowed with Hausdorff-metric and natural equivalence relation $[A]=[B]$ if $\overline{A}=\overline{B}$}.
Furthermore, if $V_i\in\mathcal{V}_n$ such that $V_i\to V\in\mathcal V_n$ as $i\to\infty$, it clearly follows that $dim(V_i)=dim(V)$ for $i$ large enough. This implies that for each $0\leq k\leq n$ we have that
\begin{equation}\label{huom1}
\{V\in\mathcal{V}_n\,:\,dim(V)=k\}=: \mathcal{V}_n^k
\end{equation}
is a compact subspace of $\mathcal{V}_n\,$. 

In the proof of Theorem \ref{OK} the differentiability of a given function has to be considered in respect to certain type of semi-linear subspaces. Therefore, let us say that $L:W\mapsto\re$, where $W\in\mathcal{W}_n$, is linear if 
\begin{equation}\label{lin1}
L(\lambda_1 w_1+\lambda_2w_2)=\lambda_1L(w_1)+\lambda_2L(w_2)\, \text{ if } w_1,w_2\in W \,,\,\,\,\lambda_1,\lambda_2\geq 0\,.
\end{equation}
The set of these mappings is denoted by $\mathcal{L}(W)\,$. As expected, $\mathcal{L}(W)$ coincides in a suitable way with the set of $W$-restrictions of linear functions $L:\rn\to\re\,$. This is briefly verified in the following proposition: 
\begin{proposition}\label{rutto}
Let $W\in\mathcal{W}_n$ and $L:W\to\re$ linear in the sense of (\ref{lin1}). Then there exists a linear function $L':\rn\to\re$ such that $L'_{|W}=L\,$. Moreover, $L'$ can be chosen such that the Lipschitz-constant of $L'$ is equal to the Lipschitz-constant of $L\,$.  
\end{proposition}
\textit{Proof.} Let $W=\langle A\rangle^+\,$, $A\subset\rn\,$. It follows from elementary linear algebra that the interior of $\langle A\rangle^+$ is nonempty in linear subspace $\langle A\rangle=:V\,$, thus there exists $B(a,r)$, $a\in A$ and $r>0$, such that $B(a,r)\cap V=:V'\subset \langle A\rangle^+$. Then it is easy to see that (\ref{lin1}) guarantees that restriction $L_{|V'}$ has a unique linear extension $L'$ to the whole linear subspace $V$. Furthermore, (\ref{lin1}) also implies  that two $W$-linear functions coincide if they coincide in an open subset of $W$, thus we get that $L'$ is also the extension of $L$. Finally, the desired extension in whole $\rn$ is naturally given by defining for $x=x_1+x_2$, $x_1\in V$ and $x_2\in V^{\perp}\,$ that
$L'(x_1+x_2)=L'(x_1)\,$.
\hfill$\Box$\\

In the proof of Theorem \ref{OK}, the following measure of the non-differentiability has a key role: 
\begin{definition}\label{tiutau2}
Let $f:\rn\to\re$, $W\in\mathcal{W}_n$, $x\in\rn$ and define
\begin{equation}\label{tiutau}
\tau(W,f,x)=\inf_{L\in\Li(W)}\,\bigg(\limsup_{w\to 0\,,w\in \overline{W}}\frac{|f(x+w)-f(x)-L(w)|}{|w|}\,\bigg).
\end{equation}
\end{definition}
Proposition \ref{rutto} guarantees that the equivalent definition of $\tau(W,f,x)$ would be also achieved if $\Li(W)$ in (\ref{tiutau}) is replaced with $\mathcal{L}(\rn)\,$.

We say that $f$ is $W$-\textit{differentiable} at $x$ if $\tau(W,f,x)=0$. One can show by an easy compactness argument that then there exists $L\in\mathcal{L}(W)$ and $\varepsilon:W\to\re$ such that $\lim_{w\to 0}\varepsilon(w)=0$ and 
\begin{equation}
f(x+w)=f(x)+L(w)+|w|\varepsilon(w)\,.
\end{equation}
In that case we say that $L$ is the \textit{$W$-derivative} of $f$ at $x\,$.

In the following lemma we verify the basic continuity properties for the directional derivatives and the functions $\tau(W,f,x)$.
\begin{lemma}\label{lem1}
Let $f$ be a Lipschitz function with Lipschitz constant $K>0$. Then\\
\text{(1)} If $f$ is directionally differentiable at $x$, then the mapping $\theta\to D_{\theta}f(x)$ on $S^{n-1}$ is Lipschitz with constant $K$. Moreover,
\begin{equation*}
\sup_{|h|\leq r}\bigg|\frac{f(x+h)-f(x)}{|h|}-D_{h}f(x)\bigg|\to 0\,\text{ as }r\to 0\,.
\end{equation*}

\noindent\text{(2)} 
The function $W\to\tau(W,f,x)$ is Lipschitz in $\mathcal{W}_n$ with constant $5K$, for every $x\in\rn\,$.
\end{lemma}
\textit{Proof.}
(1) We may assume that $x=0$ and $f(0)=0$. For the first claim, let $\theta,\theta'\in S^{n-1}$ and $\varepsilon>0$. By choosing $r$ small enough we get that
\begin{align*}
|D_{\theta}f(0)-D_{\theta'}f(0)|\leq &\bigg|\frac{f(r\theta )}{r}-\frac{f(r\theta')}{r}\bigg|+\varepsilon\\
=&\frac{|f(r\theta)-f(r\theta')|}{r}+\varepsilon\leq K|\theta-\theta'|+\varepsilon\,,
\end{align*}
where $K$ is the Lipschitz constant of $f$. This implies the first claim. For the second claim, assume, on the countrary that there exist $x_k\in\rn$, $x_k\to 0$ as $k\to\infty$, such that
\begin{equation}\label{juju}
\bigg|\frac{f(x_k)}{|x_k|}-D_{x_k}f(0)\bigg|\geq\lambda>0\,.
\end{equation}
We may assume that $\frac{x_k}{|x_k|}=:\theta_k\to \theta_0\in S^{n-1}$. The desired contradiction then follows, since the lefthandside in (\ref{juju}) can be estimated by 
\begin{align*}
&\bigg|\frac{f(|x_k|\theta_k)}{|x_k|}-\frac{f(|x_k|\theta_0)}{|x_k|}\bigg|+
 \bigg|\frac{f(|x_k|\theta_0)}{|x_k|}-D_{x_k}f(0)\bigg|\\
\leq\,& K|\theta_k-\theta_0|+\bigg|\frac{f(|x_k|\theta_0)}{|x_k|}-D_{\theta_0}f(0)\bigg|
+|D_{\theta_0}f(0)-D_{\theta_k}f(0)|\overset{k\to\infty}{\longrightarrow} 0.
\end{align*}
\\
(2) We may assume that $x=0=f(0)\,$. In order to estimate the 
difference $|\tau(W,f,0)-\tau(W',f,0)|$ for given $W,W'\in\mathcal{W}_n$, we may assume by symmetry that $\tau(W,f,0)>\tau(W',f,0)$. Let then $\varepsilon>0$ and choose a linear mapping $L:\rn\to\re$ such that 
\begin{equation*}
\tau(W',f,0)\geq\,\,\limsup_{w'\to 0\,,w'\in W'}\frac{|f(w')-L(w')|}{|w'|}\,-\varepsilon\,.
\end{equation*}
It is clear that we may assume that the Lipschitz constant of $L$ is less or equal than $3K$. For simplicity, let us denote
\begin{equation*}
G(x):=\frac{|f(x)-L(x)|}{|x|}\,.
\end{equation*}
By applying the same linear approximation $L$ also for $W$, it follows that
\begin{align}
\notag&\tau(W,f,0)-\tau(W',f,0)\leq \limsup_{x\to 0\,,x\in W}G(x)\,\,-\limsup_{x\to 0\,,x\in W'}G(x)\,+\varepsilon\\
\label{arviovio}=&\lim_{r\to 0}\bigg(\,\sup_{x\in W,|x|\leq r}G(x)\,\,-\sup_{x\in W',|x|\leq r}G(x)\,\bigg)+\varepsilon\,.
\end{align}
Suppose that above the supremum in the lefthandside is achieved, up to error at most $\varepsilon$, with $x_r\in W$, $|x_r|\leq r\,$. Then by the definition of the Hausdorff-distance and the semilinear structure of $W$ and $W'$, we find $x'_r\in W'$ so that $|x'_r|=|x_r|$ and 
\begin{equation*}
|x_r-x_r'|\leq |x_r|2\mathcal{H}_c(W,W')\,.
\end{equation*}
Therefore,
\begin{align*}
|G(x_r)-G(x'_r)|&=\,\bigg|\,\frac{|f(x_r)-L(x_r)|}{|x_r|}-\frac{|f(x'_r)-L(x'_r)|}{|x'_r|}\,\bigg|\\
&\leq\,\frac{|f(x_r)-f(x'_r)|}{|x_r|}+\frac{|L(x_r)-L(x'_r)|}{|x_r|}\\
&\leq 5K\mathcal{H}_c(W,W')\,,
\end{align*}
where $K$ is the Lipschitz constant of $f$. Combining this with (\ref{arviovio}) (and recalling that $\varepsilon$ was arbitrary small) gives the claim. 
\hfill$\Box$\\
 
The following lemma has a crucial role in the proof of our main Theorem \ref{OK} below. It is basically a simple modification of the fact that if continuous function $f$ is linear in a given half space $V+\langle b\rangle^+=:H(V,b)=:H$ and directionally linear in respect to the fixed point $a\not\in V$, then $f$ is linear also in the larger (usually) space $H-\langle a\rangle^+\,$. The proof of this fact is just a straightforward calculation. Unfortunately, unless the following Lemma is essentially based on this calculation, the proof is rather long and tedious, thanks to certain error terms, which have to be carried through the calculations. 

Another reason for the length of the following proof is that cases $a\in \langle V,b\rangle$ and $a\not\in \langle V,b\rangle$ are kind of  different nature. Remark that if $a\in \langle V,b\rangle\setminus H$, then $H-\langle a\rangle^+=H$, thus the claim is trivial. In turn, if $a\in H$ then $H-\langle a\rangle^+=\langle V,b\rangle$ and thus the linearity is extended from half-space $H(V,b)$ to linear subspace $\langle V,b\rangle\,$. However, in this case the dimension of the set of linearity is not increased, unlike in the case $a\not\in \langle V,b\rangle\,$. 
 
\begin{lemma}\label{nice}
Suppose that $f$ is Lipschitz, directionally differentiable at $x_0$ and $H(V,b)$-differentiable, where $V\in\mathcal{V}_n$, $1\leq dim(V)\leq n-1\,$, $b\in V^{\perp}$. Let $x_i\to x_0$ so that $\frac{x_i-x_0}{|x_i-x_0|}\to\theta^0\not\in V$ and 
\begin{equation}\label{olet1}
\bigg|\frac{f(x_i+h)-f(x_i)}{|h|}-D_{h}f(x_i)\bigg|<\delta
\end{equation}
for every $i$ and $h<r_0\,$. Then $\tau(H(V,b)-\langle\theta^0\rangle^+,f,x_0)\leq 5\delta\,$.  
\end{lemma}
\textit{Proof.}
We may assume that $x_0=0$ and $f(x_0)=f(0)=0$. Denote by $L_H$ the $H$-derivative of $f$ at $0$ and let $W_i=H-\langle\theta_i\rangle^+$. Moreover, let $x_i=r_i\theta_i$, where $\theta_i\in S^{n-1}$, $r_i>0$ and suppose that $w\in W_i$, thus
\begin{equation}\label{trash}
w=h+\lambda\theta_i\,\text{, where }h\in H(V,b)\text{ and }\lambda\leq 0\,.
\end{equation}

Then observe, by elementary calculation that   
\begin{equation*}
x_i+\frac{r_i}{r_i-\lambda}(w-x_i)=\frac{r_i}{r_i-\lambda}h\,=:\tilde{w}\,.
\end{equation*}
Thus, $\tilde{w}\in H$ lies on the line segment between $w$ and $x_i$. 
Then we straightforwardly calculate that
\begin{align}
\notag f(w)&=f(x_i)+\frac{f(w)-f(x_i)}{|w-x_i|}|w-x_i|\\
&=f(x_i)+\frac{f(\tilde{w})-f(x_i)}{|\tilde{w}-x_i|}|w-x_i|+R^i_1(w)\notag\\
&=f(x_i)+(f(\tilde{w})-f(x_i))\frac{r_i-\lambda}{r_i}+R^i_1(w)\notag\\
&=f(\tilde{w})(\frac{r_i-\lambda}{r_i})+\frac{\lambda}{r_i}f(x_i)+R^i_1(w)\notag\\
&=L_H(\tilde{w})(\frac{r_i-\lambda}{r_i})+\frac{\lambda}{r_i} f(r_i\theta_i)+R^i_1(w)+R_2^i(w)\notag\\
&=L_H(h)+\lambda D_{\theta^0}f(0)+R_1^i(w)+R_2^i(w)+R^i_3(w)\,,\label{nono}
\end{align}
where
\begin{align*}
R^i_1(w)&=\,\bigg|\frac{f(w)-f(x_i)}{|w-x_i|}-\frac{f(\tilde{w})-f(x_i)}{|\tilde{w}-x_i|}\bigg||w-x_i|\,,\\
R^i_2(w)&=\big(\frac{r_i-\lambda}{r_i}\big)(f(\tilde{w})-L_H(\tilde{w}))
\,\text{ and }\\
R^i_3(w)&=\lambda(\frac{f(r_i\theta_i)}{r_i}-D_{\theta^0}f(0))\,.
\end{align*}

Let us then estimate the above error terms in the case $|w|=|x_i|$. For $R_1$, notice that $|w|=|x_i|$ implies that $|\tilde{w}|\leq |w|$, since $\tilde{w}$ lies on the line segment between $w$ and $x_i$. Then it follows directly from the assumption (\ref{olet1}) that 
\begin{equation*}
|R^i_1(w)|\leq 2\delta|w-x_i|\leq 4\delta|w|\,\text{ if } w\in W_i\,,\,|w|=|x_i|\,. 
\end{equation*}
In the case of $R_2$, recall also that $|\tilde{w}|\leq |w|$. Because $L_H$ is the $H$-derivative of $f$ at $0$, it follows that (if $|w|=|x_i|$)
\begin{equation*}
\frac{|R^i_2(w)|}{|w|}\leq \big(\frac{r_i-\lambda}{r_i}\big)\varepsilon_i\,= \big(1+\frac{|\lambda|}{r_i}\big)\varepsilon_i\,,
\end{equation*}
where $\varepsilon_i\to 0$ as $i\to\infty\,$. Especially, $\varepsilon_i$ can be chosen to be independent on $w$. 
The corresponding estimate holds for $R_3$, as well. This follows, because 
\begin{align*}
\big|\frac{f(r_i\theta_i)}{r_i}-D_{\theta^0}f(0)\big|
&\leq \big|\frac{f(r_i\theta_i)}{r_i}-\frac{f(r_i\theta^0)}{r_i}\big|+\big|\frac{f(r_i\theta^0)}{r_i}-D_{\theta^0}f(0)\big|\\
&\leq K|\theta_i-\theta^0|+\big|\frac{f(r_i\theta^0)}{r_i}-D_{\theta^0}f(0)\big|\to 0\,\text{ as }i\to\infty\,.
\end{align*}

Summing up, we have shown that there exists a sequence $\varepsilon_i\overset{i\to\infty}{\longrightarrow} 0$ such that if $\lambda\leq 0$, $h\in H$ and $|h+\lambda\theta_i|=|x_i|=r_i$, then 
\begin{equation}\label{nytriitti}
\frac{|f(h+\lambda\theta_i)-L_H(h)-\lambda D_{\theta^0}f(0)|}{|h+\lambda\theta_i|}\leq 4\delta +\big(1+\frac{|\lambda|}{|h+\lambda\theta_i|}\big)\varepsilon_i\,.
\end{equation}
Furthermore, since
\begin{equation*}
\sup_{|x|\leq r\,,\,0\leq t\leq 1}\bigg|\,\frac{f(tx)}{|tx|}-\frac{f(x)}{|x|}\,\bigg|\to 0\text{ as }r\to 0\,
\end{equation*}
by Lemma \ref{lem1}, it is easy to see that (\ref{nytriitti}) holds (after a possible redefinition of $\varepsilon_i$) also if the assumption $|h+\lambda\theta_i|=|x_i|\,$ is replaced with $|h+\lambda\theta_i|\leq |x_i|\,$.
To finish the proof, we treat the cases $\theta^0\not\in \langle V,b\rangle$ and $\theta^0\in\langle V,b\rangle$ separately.

\subsection*{Case $\theta^0\not\in\langle V,b\rangle$}
Observe that assumptions $\theta^0\not\in\langle V,b\rangle$ and $\theta_i\to\theta^0$, as $i\to\infty$, guarantee that if $i$ is large enough, then each $w\in W_i$ has a unique expression $w=h+\lambda\theta_i$. Then it is easy to check that mapping $L:W_i\to\re$, defined by
\begin{equation*}
L(w)=L(h+\lambda\theta_i)=L_H(h)+\lambda D_{\theta^0}f(0)\,,
\end{equation*}
is well defined $W_i$-linear mapping. Furthermore,
assumption $\theta_i\to\theta^0\not\in \langle V,b\rangle$ as $i\to\infty$ guarantees the existence of $c>0$ such that if $i$ is large enough, then 
$|h+\lambda\theta_i|\geq c|\lambda|$ (for all $h\in H\,,\,\lambda\in\re$).
Combining this with (\ref{nytriitti}) implies that for $i$ large enough, 
\begin{equation}\label{otsa}
\frac{|f(w)-L(w)|}{|w|}\leq 5\delta\text{ if }w\in W_i\,.
\end{equation}
This implies that $\tau(W_i,f,0)\leq 5\delta$. Since $W_i\to H-\langle\theta^0\rangle^+$ as $i\to\infty$, the desired estimate $\tau(H-\langle\theta^0\rangle^+,f,0)\leq 5\delta$ follows from Lemma (\ref{lem1}).
\subsection*{Case $\theta^0\in\langle V,b\rangle$} 
Observe that in this case, if $\theta^0\not\in H(V,b)$, then $H(V,b)-\langle\theta^0\rangle^+=H(V,b)$ and the claim is trivial. Therefore, it suffices to consider the case $\theta^0\in H(V,b)$, whence $H(V,b)-\langle\theta^0\rangle^+=\langle V,b\rangle\,=\langle V,\theta^0\rangle=:W$ (recall $\theta^0\not\in V$). Then define $L:W\to\re$ by
\begin{equation*}
L(v+\lambda\theta^0)=L_H(v)+\lambda D_{\theta^0}f(0)\,,\text{ where }v\in V \text{ and } \lambda\in\re\,.
\end{equation*}
Now $L$ is well defined $W$-linear mapping and extends the $H(V,b)$-derivative $L_H$ to $\langle V,b\rangle\,$. For the claim $\tau(W,f,0)\leq 5\delta$, it suffices to show that
\begin{equation}\label{maali}
\limsup_{w\to 0\,,\,w\in W\setminus H}\frac{|f(w)-L(w)|}{|w|}\leq 5\delta\,.
\end{equation}

For (\ref{maali}), observe first that $\theta^0\not\in V$ and $\theta_i\to\theta^0$, as $i\to\infty$, guarantee the existence of $c'>0$ such that if $i$ is large enough, then $|v+\lambda\theta^0|$, $|v+\lambda\theta_i|\geq c'|\lambda|$ for all $v\in V\,$, $\lambda\in\re$.
Moreover, if $i$ is large enough, then  
\begin{align*}
|v+\lambda\theta_i|&\leq |v+\lambda\theta^0|+|\lambda||\theta^0-\theta_i|\leq |v+\lambda\theta^0|(1+\frac{1}{c'}|\theta^0-\theta_i|)\\
&\leq 2|v+\lambda\theta^0|\,\text{ for every }v\in V\,,\,\lambda\in\re\,.
\end{align*}

Combining the above facts with (\ref{nytriitti}) we finally get that if $i$ is large enough, then 
\begin{align*}
&|f(v+\lambda\theta^0)-L(v+\lambda\theta^0)|\\ \leq\,& |f(v+\lambda\theta^0)-f(v+\lambda\theta_i)|+|f(v+\lambda\theta_i)-L(v+\lambda\theta^0)|\\
\leq\,&K|\lambda||\theta^0-\theta_i|+4\delta+
\big(1+\frac{|\lambda|}{|v+\lambda\theta_i|}\big)\varepsilon_i\\
\leq\, & 5\delta|v+\lambda\theta^0|
\end{align*}
for every $v\in V$ and $\lambda<0$ such that $|v+\lambda\theta^0|\leq \frac{r_i}{2}\,$. This verifies (\ref{maali}) and completes the proof.\hfill$\Box$\\

\begin{definition}
Let $f:\rn\to\re$ and $x\in\rn$. The \textit{maximal differentiability degree} of $f$ at $x$ is defined by
\begin{equation*}
\gamma(f,x)=\max\big{\{}dim(V)\,:\,V\in\mathcal{V}_n\,,\,\tau(V+\langle b\rangle^+,f,x)=0\,\text{ for all }b\in\rn\,\big{\}}\,.
\end{equation*}
\end{definition}

\begin{theorem}\label{OK}
Let $f:\rn\to\re$ be directionally differentiable Lipschitz-function and $0\leq k\leq n-1\,$. Then 
it holds that the set
\begin{equation}
 E_k:=\{x\in\rn\,:\,\gamma(f,x)=k\}
\end{equation}
 is $\sigma$-$k$-tangential.
\end{theorem}
\textit{Proof.} Fix $0\leq k\leq n-1\,$ and let $x\in E_k$. Observe first that $\gamma(f,x)=k$ implies that there exists $\delta(x)>0$ such that for every $k+1$-dimensional linear subspace $V$ either $\tau(V,x,f)\geq \delta(x)$(especially in the case $k=n-1$) or then one can find $b\in\rn\setminus V$ such that $\tau(\langle V,b\rangle^+,f,x)\geq\delta(x)\,$. This can be shown easily by using compactness of $\mathcal{W}_n$ and $\mathcal{V}_n^k\,$, and the continuity of the function $W\to\tau(W,f,x)$ in $\mathcal{W}_n\,$ (Lemma \ref{lem1}). Therefore, we may write 
\begin{equation*}
E_k\subset\bigcup_{j=1}^{\infty}\{x\in E_k\,:\,\delta(x)>\frac{1}{j}\,\}\,=:\bigcup_{j=1}^{\infty}E_k^j\,.
\end{equation*}
Furthermore, for every $x\in E_k^j$ there exists $r(j,x)>0$ such that 
\begin{equation*}
\big|\frac{f(x+h)-f(x)}{|h|}-D_{h}f(x)\big|<\frac{1}{10j}\,,
\end{equation*}
whenever $|h|\leq r(j,x)\,$. Thus, we have
\begin{equation*}
E_k^j=\bigcup_{i=1}^{\infty}\{x\in E_k^j\,:\,r(j,x)\geq\frac{1}{i}\}=:\bigcup_{i=1}^{\infty}E_k^{j,i}\,.
\end{equation*}
Now we have reached the sufficient level of separation, thus it turns out that each set $E_k^{j,i}$ is $k$-tangential and the desired tangential $k$-plane at $x$ is exactly the $k$-dimensional linear subspace $V_x$ for which the maximal degree of differentiability is reached. To show this, suppose, on the countrary that $x,x_l\in E^{i,j}_k$ so that $x_l\to x$ and $|V_x-\frac{x_l-x}{|x_l-x|}|\not\to 0$ as $l\to\infty$. After the possible choice of a subsequence, we may assume that \begin{equation}\label{vastaol}
\frac{x_l-x}{|x_l-x|}\to\theta\in S^{n-1}\setminus V_x \text{ as } l\to\infty\,.
\end{equation}

Let then $b\in\langle V_x,\theta\rangle^{\perp}$. Let $\varepsilon >0$ and define 
\begin{equation}
H_{\varepsilon}=H(V_x,\theta+\varepsilon b)-\langle\theta\rangle^+\,.
\end{equation}
Then recall that by the definition of $\gamma(f,x)$ (and $V_x$) it holds that $\tau(H(V_x,\theta+\varepsilon b),f,x)=0$ (also if $b=0$). This, combined with $x_l\in E^{j,i}_k$, implies that 
the assumption (\ref{olet1}) in Lemma \ref{nice} is valid with constant $\delta=\frac{1}{10j}$
and we get by Lemma \ref{nice} that 
\begin{equation*}
\tau(H_{\varepsilon},f,x)\leq 5\frac{1}{10j}=\frac{1}{2j}\,.
\end{equation*}
Now it is easy to check that $H_\varepsilon\to\langle V_x,\theta\rangle+\langle b\rangle^+$ in $\mathcal{W}_n\,$ as $\varepsilon\to 0\,$ and we obtain (by Lemma \ref{lem1}) that
\begin{equation}\label{jesjes}
\tau(\langle V_x,\theta\rangle+\langle b\rangle^+,f,x)\leq \frac{1}{2j}\,.
\end{equation} 
This applies to all $b\in\langle V_x,\theta\rangle^{\perp}$. Observe that above argument is valid also in the special case $b=0$, whence $\langle V_x,\theta\rangle+\langle b\rangle^+=\langle V_x,\theta\rangle\,$
 (this corresponds the case $\theta^0\in H(V,b)$ in Lemma \ref{nice}).
Since $dim(\langle V_x,\theta\rangle)=k+1$, (\ref{jesjes}) contradicts 
with the assumption $\delta(x)\geq \frac{1}{j}\,$. The proof is complete.\hfill$\Box$\\

\subsection{Direct applications}
It is obvious that there are many important special functions in analysis, for which Theorem \ref{OK} can be directly applied. Below some basic examples are given.

\subsubsection*{Convex functions}
The estimates for the size and structure of the set of non-differentiability of a convex function has been kind of completely employed in \cite{Z} and \cite{P}. However, since 
every convex function $f:\rn\to\re$ is obviously directionally differentiable and locally Lipschitz, Theorem \ref{OK} also applies to convex functions. Kind of surprisingly, our result, while proven for much more general class of functions, is essentially optimal also for convex functions.
 
\subsubsection*{Distance functions}
Let $A\subset\rn$ non-empty and closed set and consider the distance function $g$ on $\rn\setminus A$, 
\begin{equation}\label{distanssi}
g_A(x)=\inf_{a\in A}|x-a|=-\sup_{a\in A}(-|x-a|).
\end{equation}
It is easy to see that $g_A$ is directionally differentiable. More precisely, for every $x\in\rn\setminus A$ there exists a compact set $\R(x)\subset \partial A$ of those $a\in\partial A$ for which the infimum in (\ref{distanssi}) is achieved and, especially,
\begin{equation}\label{derivdistans}
D_{\theta}g_A(x)=\min_{y\in \R(x)}\theta\cdot \frac{x-y}{|x-y|}\,
\end{equation}
for $\theta\in S^{n-1}\,$. 

By (\ref{derivdistans}) one easily obtains that $g_A$ is differentiable at $x$ exactly if $\R(x)$ is singleton. Theorem \ref{OK} thus implies that if $x\in\re^n\setminus A$ lies outside a $\sigma$-tangential exceptional set, then $x$ has the unique closest point $a_x\in A$.     
If desired, Theorem \ref{OK} actually gives us even much more precise insight into that issue. Consider, for example, the distance function for $A\subset\re^3$ and the set of points $x\in A^c$ for which there exists at least three closest points in $A$. It is almost a direct consequence of Theorem \ref{OK} that this set has to be $\sigma$-$1$-tangential.

\subsubsection*{Infimal convolution}
Let $u:\rn\to\re$ continuous, $f\in C^1(\rn\times\rn)$ (with suitable growth conditions) and define \[g(x)=\inf_{y\in\rn}\,\big(\,u(y)+f(x,y)\,\big)\,.\]
Here $g$ is so called infimal convolution of $u$ respect to function $f$, appearing as a standard tool e.g. in the theory of viscosity solutions (see e.g. \cite{CIL}).
It is rather easy to see that $g$ is directionally differentiable Lipschitz-function, thus Theorem \ref{OK} can be applied to $g$.

\subsection{Note on the optimality of Theorem \ref{OK}}
What can we say about the sharpness of Theorem \ref{OK}? One may, for example, consider the covering of the non-differentiability points even by some 'smoother' $n-1$-dimensional sets than was obtained above. The following example suggests that one can not, in general, expect remarkable refinements in this direction.\\ 

\textbf{Example.} Let us say that $\Gamma:[a,b]\to\re^2$ is regular $C^1$-curve if $|\Gamma'(t)|>0$ for every $t\in[a,b]$ and $\frac{\Gamma'(t)}{|\Gamma'(t)|}$ is continuous.
It turns out that there exists a convex function $F:\re^2\to\re^2$ so that the set of the  non-differentiability points of $F$, denoted by $E$, can not be covered by any countable union of regular $C^1$-curves.

More precisely, define a Lipschitz-continuous convex function $f:\re\to\re$ with $f(x)=0$ on $(-\infty,0]$ so that the derivative $f'$ does not exist at any rational point $q\in(0,\infty)$. Then it is easy to construct a convex function $F$ on $\re^2$ so that $F$ is not differentiable at the graph of $f$, denoted by $\mathcal{G}_f\,$.  

 For the claim, it suffices to show that $\mathcal{G}_f$ can not be covered by any desired countable collection of curves. To verify this, choose an arbitrary collection of regular $C^1$-curves $\Gamma_i:I_i\to\re^2$ and observe first that by using the above defined properties we may assume that mappings $\Gamma_i$ are injective(by possible division of $I_i$ into countable many subintervals). Moreover, by using the density of the non-differentiability points one finds for all $\Gamma_i$ and interval $I=[a-r,a+r]\subset(0,\infty)$ an interval $I'=[b-r',b+r']\subset [a-\frac{r}{2},a+\frac{r}{2}]$ so that $\Gamma_i(I_i)\cap\{(x,f(x)):x\in I'\}=\emptyset\,$. Iterating this property one can easily find a point $(x_0,f(x_0))$ which does not lie in any $\Gamma_i(I_i)\,$.

\section{Applications to the maximal functions}\label{last}
In this section we focus, in the light of Theorem \ref{OK}, on the differentiability properties of Hardy-Littlewood maximal functions. To apply Theorem \ref{OK} in a desired way, the directional differentiability of a maximal function at $x_0$ has to be considered in the situation, where essentially nothing else is assumed for $f$ than at most the differentiability at $x_0$. The additional assumptions on the continuity of the derivative or behaviour at infinity would make the arguments shorter but lead to weaker results. Taking the advantage of the local nature of our arguments, we will e.g.\ show in this connection that maximal operator preserves the almost everywhere differentiability.

\subsection{Auxiliary results}
Before the following, kind of elementary lemmata, let us introduce some notation. For the forthcoming results, let us define the \textit{restricted} maximal operator $M_{\lambda}$, defined for $\lambda\geq 0$ simply by 
\begin{equation}\label{restricted}
M_{\lambda}f(x)=\sup_{r>\lambda}\avr{B(x,r)}|f(y)|\,dy\,.
\end{equation}
Thus, using this notation $M_0f=Mf$.
It is well known that $M_{\lambda}f$, for $\lambda >0$, is Lipschitz with constant $\frac{C_n}{\lambda}\,$ for any measurable function $f$ such that $Mf\not\equiv\infty$.

For a locally integrable function $f$, the integral average of $f$ in $B(x,r)$ will be denoted by $f_r(x)$. Moreover, for technical reasons, we make a convention $f_0(x)=f(x)$ and $|f|_{\infty}(x)=\inf_{y\in\rn}M_{\lambda}f(y)\,$. We will use the following concept (introduced in \cite{Lu}) of the 'best radii' at point $x\in\rn\,$:  
\begin{definition}
Let $f:\rn\to\re$ be locally integrable and define
\[\R f(x)=\{\lambda\leq r\leq \infty\,:\,M_{\lambda}f(x)=|f|_r(x)\,\}\,.\]
\end{definition}
It turns out that $\R f(x)$ is non-empty and closed, if $\lambda>0$, for all $x\in\rn$ and in the case of the original maximal function ($\lambda=0$) for all Lebesgue point $x\in\rn$ of $f$. The special case $r=\infty$ is treated in Lemma \ref{ret} below. For the sake of simplicity, that lemma as well as the following proposition are formulated in the case $\lambda=0$ but it is obvious that exactly the same results and proofs are valid for any restricted maximal operator $M_{\lambda}$.  
\begin{lemma}\label{ret}
Suppose that $Mf(x_0)<\infty$ and $h_k\to 0$ as $k\to\infty$ so that there exists $r_k\in\R f(x_0+h_k)$ so that $r_k\to\infty$ as $k\to\infty\,$. Then $\infty\in\R f(x_0)$, thus $x_0$ is the global minimum of $Mf$, and  
\begin{equation}\label{ehto2}
\frac{Mf(x_0+h_k)-Mf(x_0)}{|h_k|}\to 0\text{ as }k\to\infty\,.
\end{equation}
\end{lemma}
\textit{Proof.}
The claim follows by rather standard estimate for Hardy-Littlewood maximal functions:  
Let $x\in\rn$ and observe that
\begin{align*}
Mf(x)&\geq\avr{B(x,r_k+|x_0+h_k-x|)}|f(y)|\,dy\\
&\geq \bigg(\frac{r_k}{r_k+|x_0+h_k-x|}\bigg)^n\avr{B(x_0+h_k,r_k)}|f(y)|\,dy\,\\
&\geq (1-\frac{|x_0+h_k-x|}{r_k+|x_0+h_k-x|})^n(Mf(x_0+h_k)\,.\\
\end{align*}
First of all, since $Mf(x_0)<\infty$, we get that $Mf(x_0+h_k)<C<\infty$ for $k$ large enough.
Then it is easy to see that the latter claim (\ref{ehto2}) follows from above by substitution $x=x_0$. 
Combining this with the lower semicontinuity of $Mf$, we obtain that $Mf(x_0+h_k)\to Mf(x_0)$ as $k\to\infty$. This finally implies, by the above estimate that $Mf(x_0)\leq Mf(x)$ for any $x\in\rn$. This completes the proof.
\hfill$\Box$\\ 

We also need the following proposition:
\begin{proposition}\label{pakkopulla}
Suppose that $Mf\not\equiv\infty$ and $f$ is continuous at $x$. Assume also that $r_i\in\R f(x+h_i)$ such that $h_i\to 0$ and $r_i\to r_0<\infty$ as $i\to\infty$. Then $r_0\in\R f(x)$.  
\end{proposition}
\textit{Proof.} Assumption $Mf\not\equiv\infty$ and the continuity of $f$ at $x$ imply that $Mf$ is continuous at $x$. Therefore, $Mf(x+h_i)\to Mf(x)$ as $i\to\infty$. But this implies the claim, since 
\begin{equation*}
Mf(x+h_i)=\avr{B(x+h_i,r_i)}|f(y)|\,dy\,\to\avr{B(x,r_0)} |f(y)|\,dy\,\text{ as }i\to\infty\,.\,\,\,\hfill\Box
\end{equation*}

The following lemma deals with the critical case $Mf(x)=|f(x)|$ and it has an important role in theorems \ref{differentioituva} and \ref{differentioituva2}. 
\begin{lemma}\label{fk}
Let $f:\rn\to\re$ be locally integrable, $x\in\rn$, and suppose that 
\begin{equation}\label{siisti}
f(y)=f(x)+D\cdot (y-x)+|y-x|u(y-x)\,,
\end{equation}
for some $D\in\rn$ and $u:\rn\to\re\,$. Then
\begin{align*}
\bigg|\avr{B(x+h,r)}f(y)\,dy\,-\avr{B(x,r)}f(y)\,dy\,-D\cdot h\,\bigg|\leq |h|C_n(\sup_{a\leq r+|h|}|u(a)|)\,.
\end{align*}
\end{lemma}
\textit{Proof.}
Observe that
\begin{align}
\notag &\avr{B(x+h,r)}f(y)\,dy\,-\avr{B(x,r)}f(y)\,dy\,\\
\label{siisti2} =&\,\frac{1}{|B_r|}\bigg(\int_{B(x+h,r)\setminus B(x,r)}f(y)\,dy\,-\int_{B(x,r)\setminus B(x+h,r)}f(y)\,dy\,\bigg)\,.
\end{align}
Plugging in the formula (\ref{siisti}) to the above integrals and using the fact 
\begin{equation*}
\frac{1}{|B_r|}\bigg(\int_{B(x+h,r)\setminus B(x,r)}D\cdot(y-x)\,dy\,-\int_{B(x,r)\setminus B(x+h,r)}D\cdot(y-x)\,dy\,\bigg)\,=\,D\cdot h
\end{equation*}
one obtains that (\ref{siisti2}) equals with
\begin{align*}
&D\cdot h\,+\frac{1}{|B_r|}\bigg(\intop_{B(x+h,r)\setminus B(x,r)}|y-x|u(y-x)\,dy\,\\&\hspace{3cm}-\intop_{B(x,r)\setminus B(x+h,r)}|y-x|u(y-x)\,dy\,\bigg)\,,
\end{align*}
where the absolute value of the latter 'error' term is bounded from above by 
\begin{equation*}
\frac{C_n|B(x+h,r)\setminus B(x,r)|}{r^n}\big(\sup_{|a|\leq r+|h|}|u(a)|\,\big)(r+|h|)\,.
\end{equation*}
By treating separately the cases $r\leq |h|$ and $r>|h|$ one can verify that 
\begin{equation}
\frac{|B(x+h,r)\setminus B(x,r)|(r+|h|)}{r^n}\leq C'_n|h|\,.
\end{equation}
This completes the proof. \hfill$\Box$\\

\subsection{Singular set of the maximal function}\label{subsub}



Let us denote the difference quotient of $f$ at $x$, respect to $h\in\rn$, by
\begin{equation}
D^h(x):=\frac{f(x+h)-f(x)}{|h|}\,,
\end{equation}
and define the \textit{singular set} of $f$ by
\[Sf:=\bigg{\{}x\in\rn\,:\,\limsup_{h\to 0}|D^h(x)|=\infty\,\bigg{\}}\,.\]

\begin{theorem}\label{ainoa}
Suppose that $f$ is a locally integrable function so that $Mf\not\equiv\infty$. Then the singular set of $Mf$ is contained in the singular set of $f\,$.
\end{theorem}
\textit{Proof.}
Suppose that $Mf\not\equiv\infty$ and $x\in S(Mf)$, thus there exists a sequence $h_k\in\rn$, $h_k\to 0$ such that
\begin{equation}\label{oletus}
D^{h_k}(Mf)(x)=\frac{|Mf(x+h_k)-Mf(x)|}{|h_k|}\longrightarrow \infty\,\text{ if }k\to\infty\,.
\end{equation}
Let us prove the claim by contradiction, thus assume that $x\not\in S(f)$ i.e. there exists constant $C>0$ and $r_0>0$ such that $D^hf(x)<C\,\,\text{ when }|h|<r_0\,.$ 
Especially this implies that $f$ is continuous at $x$, which in turn implies that $Mf$ is continuous at $x$ if $Mf(x)<\infty$. This holds, since $Mf(x)=\infty$ would (in this case) imply that $Mf\equiv\infty\,$.
  
Then, let us choose for each $k$ radius $r_k$ which almost gives the maximum average at $x+h_k$. More precisely, choose $r_k$ such that 
\begin{equation}
Mf(x+h_k)\leq \avr{B(x+h_k,r_k)}|f(y)|\,dy\,\,+\frac{|h_k|}{k}\,.
\end{equation}
It follows that $(r_k)$ is bounded, since otherwise one can use the argument in Lemma \ref{ret} to obtain a contradiction with (\ref{oletus}). 
The same follows also in the case where $r_k>\lambda>0$ for all suitably large $k\,$. This follows e.g. by observing that in this case $Mf(x)=M_{\lambda}f(x)$ and 
\begin{equation*}
|Mf(x+h_k)-M_{\lambda}f(x+h_k)|\leq \frac{|h_k|}{k}\,.
\end{equation*}

Then it is easy to check that the Lipschitz continuity of $M_{\lambda}f$ implies the desired contradiction with assumption (\ref{oletus}). The proof is thus complete if we can reach a contradiction also in the remaining case where the the sequence $r_k$ is not bounded from below. In this case we may assume, 
by extracting a subsequence, if needed that $r_k\to 0$ as $k\to\infty$. Moreover, by the continuity of $f$ at $x$ it follows that $Mf(x)=|f(x)|$.
 
Since $\limsup_{h\to 0}D^{h}f(x)<\infty$, it holds that 
\begin{equation}\label{perus}
 f(x+h)=f(x)+u_x(h)|h|\,,
\end{equation}
where $u_x:\rn\to\re$ is bounded when $|h|<r_0\,$. 

Then we recall the estimate from Lemma \ref{fk} to obtain 
\begin{equation}
\frac{Mf(x+h_k)-Mf(x)}{|h_k|}\leq C'(n)(\sup_{B(0,r_k+h_k)}(u_x))\,+\frac{1}{k}.
\end{equation}

Even easier argument shows that also
\begin{align*}
\frac{Mf(x)-Mf(x+h_k)}{|h_k|}&=\frac{|f(x)|-Mf(x+h_k)}{|h_k|}\\
&\leq\frac{1}{|h_k|}\bigg(|f(x)|-\limsup_{r\to 0}\avr{B(x+h_k,r)}|f(y)|\,dy\,\bigg)\,\\
&\leq C'(n)(\sup_{B(0,2h_k)}(u_x))\,.
\end{align*}

Since $r_k+h_k<r_0$ if $k$ is big enough and $\sup_{B(0,r_0)}(u_x)<\infty$ it follows that $D^{h_k}(Mf)(x)\not\to\infty$. This completes the proof.
\hfill $\Box$\\


\subsection{Directional differentiability and differentiability almost everywhere}
To apply the results from Section \ref{general} we have to consider the directional differentiability of maximal functions. For that, we need the following proposition:
\begin{proposition}\label{huhuh}
If $f:\rn\to\re$ is continuous and $0<r<\infty$, then $f_r$ is $C^1$-function and 
\begin{equation}\label{derivaatta}
D_{\theta}f_r(x)=\frac{C_n}{r}\avr{\partial B(x,r)}\,f(y)\frac{\theta\cdot(y-x)}{r}\,d\mathcal{H}^{n-1}(y)\,.
\end{equation}
Moreover, the mapping $(r,x)\to Df_r(x)$ is continuous on $(0,\infty)\times\rn\,$. 
\end{proposition} 
\textit{Proof.}
The formula (\ref{derivaatta}) above is just a straightforward calculation. The latter claim follows easily from (\ref{derivaatta}). 
\hfill$\Box$\\


\begin{lemma}
Suppose that $\lambda\geq 0$ and $f:\rn\to\re$ is continuous with $M_{\lambda}f\not\equiv\infty$. Then $M_{\lambda}f$ is directionally differentiable for all $\lambda>0$, and
\begin{equation}\label{kaavava}
D_{\theta}M_{\lambda}f(x)=\sup_{r\in\R f(x)} D_{\theta}|f|_r(x)\,\text{ for every }\theta\in S^{n-1}\,.
\end{equation}
Moreover, if $f$ is differentiable at $x$, then $Mf$ is directionally differentiable at $x$ and also the formula (\ref{kaavava}) is valid at $x$ for $\lambda=0$.
\end{lemma}
\textit{Proof.}
Let $\theta\in S^{n-1}$ and observe that 
\begin{equation*}
\limsup_{h\to 0}\frac{M_{\lambda}f(x+h\theta)-M_{\lambda}f(x)}{|h|}= \lim_{i\to\infty}\frac{M_{\lambda}f(x+h_i\theta)-M_{\lambda}f(x)}{|h_i|}\,, 
\end{equation*}  
for some sequence $(h_i)$, $h_i\to 0$ as $i\to\infty$. If $\lambda>0$, $M_{\lambda}f$ is Lipschitz, thus the limit on the right hand side exists. In the case $\lambda=0$, the differentiability of $f$ at $x$ implies that $x$ does not lie in the singular set of $|f|$, thus the limit exists by virtue of Theorem \ref{ainoa}. 

Suppose then that $r_i\in\R f(x+h_i\theta)\,$. By extracting a subsequence, if needed, we may assume that $r_i\to r_0\in\R f(x)$ as $i\to\infty$, $\lambda\leq r_0\leq \infty$ (Lemma \ref{pakkopulla} is used here). Then
\begin{equation}\label{ekaeka}
\frac{Mf(x+h_i\theta)-Mf(x)}{|h_i|} 
\leq \frac{|f|_{r_i}(x+h_i\theta)-|f|_{r_i}(x)}{|h_i|}\to D_{\theta}(|f|_{r_0})(x)\,.
\end{equation}
The convergence above results from the previous auxiliary results; in the case where $0<r_0<\infty$, one can easily show that the convergence is valid by the continuity of the mapping $(r,x)\to D_{\theta}|f|_r(x)$ on $(0,\infty)\times\rn\,$, following from Proposition \ref{huhuh}. In the case of $r_0=\infty$ (recall the convention $|f|_{\infty}\equiv\inf_{y\in\rn}M_{\lambda}f(y)$) one has to use Lemma \ref{ret}. Finally, if $\lambda=0$, it may happen that $r_0=0$. In this case the convergence is valid, by Lemma \ref{fk}, if $|f|$ is differentiable at $x$. This holds if $f$ is differentiable at $x$,  since $Mf(x)=|f(x)|>0$ ($Mf(x)=0$ implies $f\equiv0$).  

For the reverse inequality, observe that for every $r\in\R f(x)$ (also if $r=0$ or $r=\infty$) it holds that 
\begin{equation}\label{tokatoka}
\liminf_{h\to 0}\frac{Mf(x+h\theta)-Mf(x)}{|h|}\geq \liminf_{h\to 0}\frac{|f|_{r}(x+h\theta)-|f|_{r}(x)}{|h|}
= D_{\theta}|f|_{r}(x)\,. 
\end{equation}
Obviously the claim follows from (\ref{ekaeka}) and (\ref{tokatoka}). Remark that in the case $\lambda >0$ one has to consider above only the case $r_0>0$ (\ref{ekaeka}) and $r\geq\lambda>0$ (\ref{tokatoka}). In this case the needed auxiliary lemmas does not assume the differentiability for $|f|$. 
\hfill$\Box$\\

Then we obtain the following corollary.
\begin{corollary}\label{koro}
If $f$ is a continuous function such that $M_{\lambda}f\not\equiv\infty$ and $\lambda>0$, then $M_{\lambda}f$ satisfies the assumptions of Theorem \ref{OK}, thus it holds for $0\leq k\leq n$ that the set where the maximal differentiability degree of $M_{\lambda}f$ equals to $k$ is at most $\sigma$-$k$-tangential. Moreover, the same conclusion holds for $Mf$ (case $\lambda=0$) if $f$ is differentiable and Lipschitz. 
\end{corollary}

\begin{lemma}\label{uusiuusi}
Let $f:\rn\to\re$ such that $Mf\not\equiv\infty$ and let $E_{|f|}$, $E_{j}$ and $E_{Mf}$ denote the non-differentiability points of $|f|$, $\max\{|f|,M_{\frac{1}{j}}f\}$ and $Mf$, respectively. Then
\begin{equation}
 E_{Mf}\subset E_{|f|}\cup\bigcup_{j=1}^{\infty}E_j\,.
\end{equation} 
\end{lemma}
\textit{Proof.}
Suppose that $x\in E_{Mf}\setminus E_{|f|}$. Since $x\not\in E_{|f|}$, $|f|$ is continuous(even differentiable) at $x$. Therefore, if $Mf(x)>|f(x)|$, there exists $j_0\in\na$ such that $Mf$ coincides with $M_{\frac{1}{j_0}}f$ in a neigbourhood of $x$. In this case it clearly holds that $x\in E_{j_0}$. On the other hand, if $Mf(x)=|f(x)|$, then it is easy to check that $x\in E_{Mf}\setminus E_{|f|}$ implies that there exists a sequence $h_k\to 0$ as $k\to\infty$ and $c>0$ such that 
\begin{equation}\label{eiei}
 Mf(x+h_k)\geq |f(x)|+D|f|(x)\cdot h_k+c|h_k|
\end{equation}
for all $k\in\na$. Let then $r_k\in\R f(x+h_k)$, thus $Mf(x+h_k)=|f|_{r_k}(x)$. Then it holds that $\liminf_{k\to\infty}r_k > 0$, since the opposite claim would yield a contradiction with (\ref{eiei}) by virtue of Lemma \ref{fk}. This in turn implies that there exists $j_0\in\na$ such that $M_{\frac{1}{j_0}}f(x+h_k)=Mf(x+h_k)$ for $k$ large enough. Therefore, it clearly follows that $\max\{|f|,M_{\frac{1}{j_0}}f\}$ is not differentiable at $x$, thus $x\in E_{j_0}$. This completes the proof.
\hfill$\Box$\\ 

As a corrollary, we obtain that maximal operator preserves the a.e. differentiability:
\begin{theorem}\label{differentioituva2}
If $f$ is a.e. differentiable and $Mf\not\equiv\infty$, then $Mf$ is a.e. differentiable.
\end{theorem}
\textit{Proof.}
It is elementary fact that if $f$ and $g$ are differentiable a.e., then the same holds for $\max\{f,g\}$ as well. Combining this with Lemma \ref{uusiuusi} gives the claim.\hfill$\Box$\\

\textit{Remark.} Theorem \ref{differentioituva} also follows directly from Theorem \ref{ainoa} (without Lemma \ref{uusiuusi}) by using Stepanov's Theorem (\cite[3.1.8]{F}), which says that any measurable function $f:\rn\to\re$ is differentiable a.e. outside the singular set $Sf$.

\subsection{Proof of Theorem \ref{differentioituva}}
Theorem \ref{differentioituva} deals with the case where the Lipschitz-assumption for $f$ in the latter statement of Corollary \ref{koro} is dropped, indeed the assumptions in Theorem \ref{differentioituva} were that $f:\rn\to\re$ is continuous and differentiable outside a $\sigma$-tangential set and $Mf\not\equiv\infty$. We have to show that then $Mf$ is also continuous and differentiable outside a $\sigma$-tangential set. The claim turns out to follow easily from Lemma \ref{uusiuusi} and the following elementary proposition:
\begin{proposition}\label{HUH}
 Let $f$ and $g$ be differentiable outside $\sigma$-tangential sets $E_g$ and $E_f$ (respectively). Then $\max\{f,g\}$ is differentiable outside a $\sigma$-tangential set. 
\end{proposition}
\textit{Proof.} Let $E$ denote the set where $\max\{f,g\}$ is not differentiable. Observe first that
$E\cap\{f(x)\not= g(x)\}\subset E_f\cup E_g$, and $E\cap\{f(x)=g(x)\}\cap(E_f\cup E_g)\subset E_f\cup E_g$. Therefore, it suffices to show that 
\begin{equation*}
E\cap\{f(x)=g(x)\}\setminus(E_f\cup E_g)\,\text{ is $\sigma$-tangential}\,.
\end{equation*}
This follows easily by observing that if $x$ lies in the above set, then $Df(x)$ and $Dg(x)$ exist and $Df(x)\not= Dg(x)$, implying that $\{f=g\}$ is tangential at $x\,$.
\hfill$\Box$\\

As in Lemma \ref{uusiuusi}, denote by $E_{|f|}$, $E_{j}$ and $E_{Mf}$ the non-differentiability points of $|f|$, $\max\{|f|,M_{\frac{1}{j}}f\}$ and $Mf$, respectively. Since $f$ is differentiable up to $\sigma$-tangential set, the same applies to $|f|$, as well. Thus $E_{|f|}$ is $\sigma$-tangential. Moreover, Corollary \ref{koro} guarantees that $M_{\frac{1}{j}}f$ is differential up to a $\sigma$-tangential set for each $j\in\na$, whence Proposition \ref{HUH} above tells that $E_{j}$ is $\sigma$-tangential for all $j\in\na$. Since the union of all these exceptional sets is again $\sigma$-tangential, the claim follows from Lemma \ref{uusiuusi}.
\hfill$\Box$\\ 

\subsection{Other maximal operators and general pointwise maximum function}
It is clear that directional differentiability holds for various other maximal functions, as well, and thus Theorem \ref{OK} is in our use, if only certain sufficient Lipschitz-conditions are satisfied. This applies, for example in the case, where the balls in the definition of $Mf$ (or $M_{\lambda}f$) are replaced with $n$-dimensional cubes, in the case of \textit{non-centered} maximal operator (where balls $B(x,r)$ in (\ref{kaava:1}) are replaced with all balls containing point $x$) or in the case of so called \textit{fractional} maximal operator.

One may also consider the following more general class of maximal type functions satisfying the assumptions of Theorem \ref{OK}: suppose that $\{f_k\}$ is a countable family of $C^1$-functions $f_k:\rn\to\re$ and define their pointwise maximum function $F$ by
\begin{equation}\label{suppi}
F(x)=\sup_{k\in\na}f_k(x)\,.
\end{equation}
This kind of functions are sometimes called in literature as \textit{regular upper envelopes} or  \textit{pointwise maximum-functions}, see e.g. \cite[Chapter 4]{BC} and references therein.
It turns out that assuming $\{Df_k\}_{k\in\na}$ to be locally uniformly bounded and equicontinuous implies that $F$ is directionally differentiable and Lipschitz.
The proof of this fact turns out to be rather elementary and it is left to the interested reader.

\end{document}